\documentclass{article}
\usepackage{amsfonts}
\usepackage{amsmath}
\usepackage{amssymb}
\newtheorem{lemma}{Lemma}[section]
\newtheorem{definition}{Definition}
\newtheorem{proposition}[lemma]{Proposition}
\newtheorem{corollary}[lemma]{Corollary}
\newtheorem{theorem}[lemma]{Theorem}

\newcommand{\udim}{\overline{\dim}}
\newcommand{\ldim}{\underline{\dim}}

\newcommand{\diam}{\textrm{diam}}
\newcommand{\supp}{\textrm{supp}}

\newcommand{\ma}{\mathcal{A}}
\newcommand{\bbn}{\mathbb{N}}

\title{A modified multifractal formalism for a class of self-similar measures with overlap}
\author{Pablo Shmerkin}
\begin{document}
\maketitle

\textbf{Abstract}.

The multifractal spectrum of a Borel measure $\mu$ in
$\mathbb{R}^n$ is defined as
\[
f_\mu(\alpha) = \dim_H\left\{x:\lim_{r\rightarrow 0} \frac{\log
\mu(B(x,r))}{\log r}=\alpha\right\}.
\]
For self-similar measures under the open set condition the
behavior of this and related functions is well-understood
(\cite{mauldin},\cite{olsen},\cite{randommultifractals}); the
situation turns out to be very regular and is governed by the
so-called ``multifractal formalism''. Recently there has been a
lot of interest in understanding how much of the theory carries
over to the overlapping case; however, much less is known in this case and what is known makes it clear that more complicated phenomena are possible. Here we carry out a complete study of the
multifractal structure for a class of self-similar measures with
overlap which includes the $3$-fold convolution of the Cantor
measure. Among other things, we prove that the multifractal
formalism fails for many of these measures, but it holds when
taking a suitable restriction.

\section{Introduction and statement of results}

Recall that given a Borel measure $\mu$ in $\mathbb{R}^n$, the
upper and lower local dimensions of $\mu$ are defined as
\[
\udim\mu(x) = \limsup_{r\rightarrow 0^+} \frac{\log \mu
(B(x,r))}{\log r};
\]
\[
\ldim\mu(x) = \liminf_{r\rightarrow 0^+} \frac{\log \mu
(B(x,r))}{\log r}.
\]
When $\udim\mu(x)=\ldim\mu(x)$ we refer to the common value as
the local dimension of $\mu$ at $x$, and we denote it by
$\dim\mu(x)$. One of the main objectives of multifractal analysis
is to study the level sets of the (upper or lower) local
dimensions of a given measure. To this end a number of
``multifractal spectra'' are introduced; perhaps the most widely
used is
\[
f_H(\alpha) = f_{\mu,H}(\alpha) = \dim_H \{x\in\supp\mu:
\dim\mu(x)= \alpha\},
\]
where $\dim_H$ is the Hausdorff dimension. In the context of multifractal analysis it is convenient to use the convention $\dim_H(\varnothing)=-\infty$, and we will implicitly do so in the sequel.

Another basic function
in multifractal theory is the $L^q$-spectrum, which is defined
(for compactly supported measures) as follows: let
\[
S_{\mu,r}(q) = \sup \left \{\sum_i \mu(B(x_i,r))^q:
\{B(x_i,r)\}_{i} \textrm{ is a packing of } \supp\mu \right\}.
\]
The (lower) $L^q$-spectrum is then given by
\[
\underline{\tau}(q) = \underline{\tau}_{\mu}(q) =
\liminf_{r\rightarrow 0} \frac{\log S_{\mu,r}(q) }{\log r}.
\]
When the limit above exists for all $q$ we speak of ``the''
$L^q$-spectrum $\tau(q)$. The celebrated heuristic principle
known as the ``multifractal formalism'' states that $f_H(\alpha)$
and $\tau(q)$ (or $\underline{\tau}(q)$) form a
Legendre-transform pair. Recall that given a concave function
$g(q):\mathbb{R}\rightarrow [-\infty,\infty)$, its Legendre
transform is given by
\[
g^*(\alpha) = \inf_{q\in\mathbb{R}} q\alpha-g(q).
\]
It is well-known that $g^*$ is also concave, and moreover
$g^{**}=g$. For different accounts and generalizations of the
multifractal formalism see \cite{mauldin}, \cite{olsen} and
\cite{pesin}.

Although false in general, the multifractal formalism has been
verified for many natural measures, including self-similar
measures under the open set condition (\cite{mauldin},
\cite{randommultifractals}). Without separation, however, much
less is known, and most of that is known refers to the portion
of the $L^q$-spectrum corresponding to $q\ge 0$ (or,
equivalently, the portion of $f(\alpha)$ corresponding to
$\alpha\le \gamma=\max f)$; see \cite{wsc} and \cite{lqspectrum}
for some of the deep results obtained. In \cite{convolution} a
first investigation was made for the case $q<0$; there the authors
proved the striking result that for the $m$-fold convolution of
the Cantor measure, $m\ge 3$, the maximum of the set of local
dimensions is an isolated point. More recently, D.J. Feng and E. Olivier \cite{weakgibbs} proved that the multifractal formalism holds under a so-called ``weak-Gibbs'' condition. In particular, they showed that this condition is verified for the Bernoulli convolution associated with the Golden ratio.

In this paper we undertake a detailed study of a family of
self-similar measures with overlap that includes the $3$-fold
convolution of the standard Cantor measure. Our main result is
that the multifractal formalism fails for these measures, but
it is verified if we restrict the measure to any subinterval of
the support not containing its extreme points, see Theorems
\ref{theo:lqspectrumpart} and \ref{theo:multifractalformalism}.
We also find formulae for the relevant extreme dimensions and
prove that the local dimension exists and is almost everywhere
constant; some of these expressions are easy to estimate numerically, while others seem to require very heavy computing power. Finally, we investigate some concrete examples which
exhibit different multifractal phenomena and suggest future
research.

From now on fix integers $d\ge 3$ (the ``base'') and $m\ge d$
(the ``maximum digit''), together with a probability vector
$\mathbf{p}=(p_0,\ldots,p_m)$. Two related numbers will appear
frequently, so we will denote them by special symbols: let
\begin{equation} \label{eq:defdelta}
\delta = -1/\log d,
\end{equation}
and
\begin{equation} \label{eq:defxi}
\xi = m/(d-1).
\end{equation}
For technical reasons we will need to impose the following
condition on $\mathbf{p}$.
\begin{definition}
The probability vector $\mathbf{p}$ is said to be \em{regular} if
\[
p_0,p_m \le p_i \textrm{ for all } i=1,\ldots, m-1.
\]
\end{definition}
We will always assume that $\mathbf{p}$ is regular. We remark,
however, that most of the results hold under weaker conditions on
$\mathbf{p}$, and some results are valid for \em every \em
probability vector $\mathbf{p}$.

Let $\mu$ be the self-similar measure associated to the weighted
iterated function system
\[
\left\{ \left(\frac{x}{d}+\frac{i}{d},p_i\right):0\le i\le m
\right\}.
\]
In other words, $\mu$ is the only compactly supported Borel
probability measure on $\mathbb{R}$ such that  $\mu(A) =
\sum_{i=0}^m p_i \,\mu(dA-i)$ for every Borel set $A$. Another
convenient way to look at $\mu$ is as the distribution of the
random sums $\sum_{j=0}^\infty X_i d^{-i}$, where $X_i$ takes the
value $i$ with probability $p_i$. It is not difficult to see that
the $k$-fold convolution of the standard Cantor measure can be
represented in this way, by taking $d=3$, $m=k$ and
$p_i=2^{-k}{k\choose i}$, see Lemma \ref{lemma:convolution}.

We will now introduce some notation. For any finite set
$\mathcal{A}$ define
\[
\begin{array}{lcl}
\mathcal{A}^k & = & \{\sigma=(\sigma_1,\ldots,\sigma_k):\sigma_j\in\mathcal{A}\};\\
\mathcal{A}^* & = & \bigcup_{k=0}^\infty \mathcal{A}^k; \\
\mathcal{A}^{\mathbb{N}} & = & \{\sigma= (\sigma_1,\sigma_2
,\ldots):\sigma_j
\in \mathcal{A} \}; \\
\overline{\mathcal{A}^*} & = & \mathcal{A}^* \cup
\mathcal{A}^{\mathbb{N}}; \\
|\sigma| & = & \textrm{ length of } \sigma\in \overline{\mathcal{A}^*} ;\\
\sigma|k & = & \textrm{ restriction of } \sigma\in
\overline{\mathcal{A}^*} \textrm{ to its first }
k \textrm{ coordinates }  ;\\
T & = & \textrm{ shift operator on } \mathcal{A}^{\mathbb{N}};\\
 \nu_{\mathbf{p}} & = & \textrm{ Bernoulli (product)
measure corresponding to } \mathbf{p} = \{p_i\}_{i\in\mathcal{A}}.
\end{array}
\]
If $\sigma^1,\ldots,\sigma^j\in\mathcal{A}^*$ we will denote their
juxtaposition by $(\sigma^1,\ldots,\sigma^j)$. We specialize now
to the case $\mathcal{A}=\{0,\ldots,m\}$, and define the
``projection'' $\pi: \ma^\bbn\rightarrow\mathbb{R}$ by
\[
\pi(\omega) = \sum_{i=1}^\infty \omega_i \, d^{-i}.
\]
A standard fact is that $\mu$ is the projection of
$\nu=\nu_{\mathbf{p}}$, in the sense that
$\mu(A)=\nu(\pi^{-1}(A))$. We will define $\pi$ also for
$\sigma\in\mathcal{A}^*$; i.e, $\pi(\sigma)=\sum_{i=1}^{|\sigma|}
\sigma_i d^{-i}$. Observe that $\supp\mu = \pi(\ma^\bbn) =
[0,\xi]$, and
\begin{equation} \label{eq:projectionatoms}
\pi(\mathcal{A}^k) = \left\{j d^{-k}:0\le j\le \sum_{i=1}^k m
d^i\right\}.
\end{equation}
If $s\in\mathbb{R}$ let
\[
[s]_k = \{ \sigma\in\ma^k:\pi(\sigma)=s \}.
\]
Note that $[s]_k$ is empty unless $s\in\pi(\mathcal{A}^k)$. For
$\sigma\in\mathcal{A}^k$ let
\[
\mathbf{p}(\sigma) = \prod_{j=1}^k p_{\sigma(j)} =
\nu(\{\omega\in\mathcal{A}^{\mathbb{N}}:\omega|k=\sigma\}),
\]
and
\[
\eta(\sigma) = \sum\{ \mathbf{p}(\sigma'): \sigma'\in
[\pi(\sigma)]_k\} = \sum\{ \mathbf{p}(\sigma'):|\sigma'|=|\sigma|, \pi(\sigma')=\pi(\sigma) \}.
\]
The function $\eta$ will play a key role in what follows. To
motivate this, observe that if we define $\mu_0=\delta_0$ (the
Delta measure giving full mass to $0$), and
\[
\mu_{k+1}(A) = \sum_{i=0}^m p_i \mu_k(d A-i),
\]
then $\mu_k$ converges weakly to $\mu$ and it can be easily
verified that $\mu_k$ assigns $\eta(\sigma)$ mass to
$\pi(\sigma)$. Actually, more is true: $\eta(\sigma)$ can be
roughly compared to $\mu(B(\pi(\sigma),d^{-|\sigma|}))$; see
Proposition \ref{prop:dimensionformula} below.

Recall that the minimum of the $p_i$ is attained at $p_0$ or
$p_m$ (or both). Notice that if we replace $\mathbf{p}$ by
$\mathbf{p}^* = (p_m,\ldots,p_0)$, then the resulting measure
$\mu^*$ is just a flipped version of $\mu$: $\mu^*(A) =
\mu(\xi-A)$ for every $A\subset\mathbb{R}$; thus $\mu$ and $\mu^*$
are indistinguishable from the multifractal point of view. Hence
we will henceforth assume, without loss of generality, that
$p_0\le p_m\le p_i$ for $i=1,\ldots,m-1$. Iterating the defining
equation $\mu(A) = \sum_{i=0}^m p_i \,\mu(dA-i)$ we obtain  that
$\mu$ is the only probability measure satisfying
\[
\mu(A) = \sum_{\sigma\in\mathcal{A}^k} \mathbf{p}(\sigma)
\mu(d^k(A-\pi(\sigma))) = \sum_{j=0}^{m+\ldots+md^{k-1}}
\eta(\sigma^j)\mu(d^k A-j),
\]
where $\sigma^j$ is a representative of $[jd^{-k}]_k$. Hence if we
replace $d$ by $d^k$, $m$ by $m(1+\ldots+d^{k-1})$ and
$\mathbf{p}$ by $\{\eta(\sigma^j)\}_j$ we get a weighted IFS with
the same attractor $\mu$ (it follows easily by induction that the
new weights are still regular if $\mathbf{p}$ is regular). By
choosing $k=2$ if necessary we can assume that $m>d$; we will do
so unless otherwise stated.

We will now introduce a set of transfer matrices. Write
\[
a = 1+ \left\lfloor \frac{m-d}{d-1} \right\rfloor.
\]
Let us adopt the convention that $p_i=0$ for
$i\notin\{0,\ldots,m\}$, and define functions
$M_0,\ldots,M_m:\{-a,\ldots,a\}^2\rightarrow\mathbb{R}$ by
\[
M_i(k,l) = p_{-ld+k+i}.
\]
These functions can be of course considered as matrices in
$\mathbb{R}^{(2a+1)\times (2a+1)}$. For $\sigma\in\mathcal{A}^k$ we
will write $M(\sigma)=M_{\sigma(k)}\cdots M_{\sigma(1)}$ (note
that $\sigma$ is reversed in the product). The use of
appropriately defined transition matrices seems to be a recurrent
tool in the investigation of self-similar measures with overlap;
see the survey \cite{wscsurvey} for some instances of this.

We recall some well-known facts and definitions from Linear
Algebra. A matrix norm $\|\cdot\|$ is \textit{consistent} if
$\|AB\|\le \|A\|\|B\|$; all operator norms are consistent. For
all consistent norms it is verified that
\begin{equation} \label{eq:radiusnorm}
\rho(A) = \lim_{k\rightarrow\infty} \|A^k\|^{1/k} \le \|A\|,
\end{equation}
where $\rho(A)$ is the spectral radius of $A$. The generalized
spectral radius of a family of matrices $\mathcal{M}$ is defined
as
\[
\tilde{\rho}(\mathcal{M}) = \limsup_{k\rightarrow\infty}
\left(\sup \{\rho(A_1 A_2\ldots A_k): A_i\in\mathcal{M}
\}\right)^{1/k};
\]
see \cite{finiteness} for a discussion of this and related
concepts. Our first result deals with the extreme local
dimensions of $\mu$. Let $\overline{\Delta}\mu = \{\udim
\mu(x):x \in \supp \mu \}$, and define $\underline{\Delta}\mu$
analogously. Let also
\[
\Delta\mu = \{\alpha:\alpha=\dim\mu(x) \textrm{ for some }
x\in\supp\mu\}.
\]
Note that $\Delta\mu\subset \overline{\Delta}\mu \cap
\underline{\Delta}\mu$, but there is no \textit{a priori} reason
for equality.

\begin{theorem} \label{theo:extremedimensions}
Let
\[
\begin{array}{lll}
\underline{\alpha} & = & \delta \log \tilde{\rho} (M_0,\ldots,M_m) ; \\
\alpha^* & = & \delta\log\underline{\rho} ({M_0,\ldots,M_m});\\
\overline{\alpha} & = & \delta \log p_0;
\end{array}
\]
where
\[
\underline{\rho}({M_0,\ldots,M_m}) = \inf\left\{ \rho(M(\sigma))^{1/|\sigma|}:\sigma\in\mathcal{A}^*\textrm{ and }
\sigma_1\notin \{0,m\}\right\};
\]
Then:
\begin{enumerate}
\item
$\inf\underline{\Delta}\mu=\inf\Delta\mu=\underline{\alpha}$;
\item
$\sup\overline{\Delta}\mu=\sup\Delta\mu=\overline{\alpha}$, and
this common supremum is attained;
\item
$\alpha^*$ is the minimum number such that
\begin{equation} \label{eq:estimatedeltamu}
\Delta\mu \subset [\underline{\alpha},\alpha^*] \cup \{\delta\log
p_m \} \cup \{\overline{\alpha}\}.
\end{equation}
\item
$\alpha^*<\overline{\alpha}$ if $p_0<p_i$ for all
$i=1,\ldots,m-1$.
\end{enumerate}
\end{theorem}

We remark that the formula for $\overline{\alpha}$ and the
interesting fact that $\overline{\alpha}$ is isolated if $p_0$ is
a strict minimum of the weights are straightforward
generalizations of results in \cite{convolution}. On the other
hand, the expressions for $\underline{\alpha}$ and $\alpha^*$ are
new results even for the $m$-fold convolution of the Cantor
measure if $m\ge 5$.

We will often need to work with subsets $\Xi_{k}$ of $\ma^k$ such
that the restriction $\pi|_{\Xi_k}$ is injective, but still
$\pi(\Xi_k)= \pi(\ma^k)$; this is equivalent to choosing one
representative from each nonempty class $[j d^{-k}]_k$. We will
henceforth assume that such a family $\{\Xi_k\}_{k\in\bbn}$ has
been selected.

\begin{proposition} \label{prop:lqspectrum}
The $L^q$ spectrum $\tau(q)$ exists for all $q$. Moreover, if
\[
\overline{S}_k(q) =  \sum_{\sigma\in\Xi_k}
\eta(\sigma)^q,
\]
then $\tau(q) = \lim_{k\rightarrow\infty} \delta k^{-1} \log
\overline{S}_k(q)$.
\end{proposition}

We impose now the additional condition $m<2d-2$ (or,
alternatively, $\supp\mu \subset [0,2)$). In this setting sharper
and more complete results can be obtained; this is due to the
availability of a ``barrier digit'', which is defined in the next
lemma.
\begin{lemma}
Assume that $m<2d-2$. There exists $\sigma\in\mathcal{A}^*$ such
that
\[
\pi(\omega)-1 < 0 < \xi < \pi(\omega) + 1
\]
for every $\omega\in\mathcal{A}^{\mathbb{N}}$ such that
$\omega||\sigma|=\sigma$. If this happens for some $\sigma=(b)$ (a
single digit), we will call $b$ a \em barrier digit \em.
\end{lemma}
\textit{Proof}. Choose $\varepsilon>0$ and $x=\pi(\omega)
\in\supp\mu$ such that
\[
\xi - (1+\varepsilon) < x < 1 - \varepsilon.
\]
Take $k$ such that $\sum_{j=k+1}^\infty m d^{-j}<\varepsilon/2$;
then $\sigma=\omega|k$ verifies the desired property.
$\blacksquare$

Note that if $m<2d-2$ then for a suitable iteration of the IFS
there is a barrier digit, so we will always assume that there is
such a digit already in the original IFS. We work with ``barrier digits''
instead of ``barrier words'' just for notational convenience. In particular, we stress that
although there is no barrier digit for the $3$-fold convolution of the Cantor measure, there is one for
some suitable iteration, so our results do apply to this important example.

The importance of the barrier digit lies in that it allows to restore a weak form of
uniqueness in the representation $\pi(\omega)=x$; see Lemma
\ref{lemma:barrierdigit} in the next section. Suppose that there
is a barrier digit $b$, and define
\begin{equation} \label{eq:shat}
\hat{S}_k(q) = \sum_{\sigma\in\Xi_{k-1}} \eta(b,\sigma)^q.
\end{equation}

($(b,\sigma)$ denotes the concatenation of the digit $b$ and the word $\sigma$). Observe that $\hat{S}_k$ does not depend on the choice of $\Xi_k$,
since $\pi(\sigma)=\pi(\sigma')$ implies that
$\eta(b,\sigma)=\eta(b,\sigma')$. We state now our main results.

\begin{theorem} \label{theo:deltamupart}
\[
(\underline {\alpha}, \alpha^*) \cup \{\delta\log p_m\}  \cup
\{\overline{\alpha}\} \subset \Delta\mu.
\]
\end{theorem}

Note that from item 3. of Theorem \ref{theo:extremedimensions} it
follows that the inclusion in this theorem is in fact an
equality, with the possible exceptions of $\underline{\alpha}$
and $\alpha^*$, which a priori may or may not be in
$\Delta{\mu}$.

\begin{theorem} \label{theo:lqspectrumpart}
The following limit exists for all $q$:
\begin{equation} \label{eq:definitionhattau}
\hat{\tau}(q) = \lim_{k\rightarrow\infty} \delta k^{-1} \log
\hat{S}_k(q).
\end{equation}
Moreover, the limit can be replaced by supremum.

If $K\subset (0,\xi)$ is any closed interval then
$\hat{\tau}$ is the $L^q$-spectrum of $\mu|_K$, the restriction
of $\mu$ to $K$.

Finally, if $\alpha^*=\overline{\alpha}$ then $\tau(q)
=\hat{\tau}(q)$ for all $q$, while if
$\alpha^*<\overline{\alpha}$, then there is $q_0\in (-\infty,0)$
such that
\begin{equation} \label{eq:formulalqspectrum}
\tau(q) = \left\{
\begin{array}{lll}
\overline{\alpha} q & \textrm{ if } & q \le q_0 \\
\hat{\tau}(q) & \textrm{ if } & q_0 < q
\end{array}
\right..
\end{equation}
\end{theorem}

This theorem is the key to the understanding of the
multifractal formalism for $\mu$: since the measure has very low
concentration near the endpoints of the support, the contribution
of one single ball centered at $0$ in the sums $S_r(q)= \sum_i
\mu(B(x_i,r))^q$  is greater than all the others, provided $q$ is
sufficiently close to $-\infty$. This instability precludes the
multifractal formalism from holding near $\alpha^*$, as $\tau(q)$
does not accurately reflect the distribution of the measure near
$q=-\infty$. Restricting the measure to a subinterval removes the
instability but, thanks to self-similarity, does not alter the
``correct'' value of the $L^q$-spectrum. Hence it is not
surprising that $f_H(\alpha)$ can be computed as the Legendre
transform of $\hat{\tau}(q)$, with the obvious exceptions of the
dimensions attained at the extreme points $0,\xi$; this is the
content of our next theorem.

\begin{theorem} \label{theo:multifractalformalism}
$f_H(\alpha) = \hat{\tau}^*(\alpha)$ for every $\alpha\in
(\underline{\alpha},\alpha^*)$, where  ${\hat{\tau}}^*$ denotes
the Legendre transform of $\hat{\tau}$.
\end{theorem}

An important feature of a multifractal measure which is not
implicit in either the multifractal or $L^q$ spectra is the
existence (or lack thereof) of an almost sure local dimension.
This question is answered in the next proposition.

\begin{proposition} \label{prop:almostsuredimension}
For any $m\ge d$ the local dimension exists and is constant
almost everywhere. The almost sure value is given by
\[
\gamma = \inf_{k}\, \frac{1}{k} \sum_{\sigma\in\Xi_k}
-\eta(\sigma)\log\eta(\sigma) = \sup_{k}\, \frac{1}{k} \sum_{\sigma\in\mathcal{A}^k} -p(\sigma) \log\|M(\sigma)\|.
\]
where $\|\cdot\|$ is a fixed consistent norm.
\end{proposition}

We remark that when an almost sure local dimension $\gamma$ exists, many other
dimensions of the measure, such as Hausdorff or entropy
dimensions, are also equal to $\gamma$.

After a first version of this paper was completed we were informed that Ka-Sing Lau and Xiang-Yang Wang obtained similar results regarding the $L^q$ spectrum, albeit using different methods. They give a formula for the $L^q$-dimension of the $3$-fold convolution of the Cantor measure (the formula is different from ours); it follows from their formula that the $L^q$-dimension is analytic except at one point. On the other hand, they do not consider some of the questions studied here (like existence of an almost sure local dimension). Their results extend to other overlapping self-similar measures, but in a different direction than the class studied here. We remark that they use the Renewal Theorem and some detailed combinatoric estimates, while our techniques have a more linear algebraic flavor. They also announce results similar to ours regarding the multifractal spectrum, but we have not been able to see their proofs yet.

\section{Auxiliary results}

This section contains the main technical ingredients of the paper.
We begin with a lemma that, although very simple, will play a
fundamental role in the sequel.

\begin{lemma} \label{lemma:supermultiplicativity}
If $\sigma,\sigma'\in\mathcal{A}^*$ then
\[
\eta(\sigma,\sigma') \ge \eta(\sigma)\eta(\sigma').
\]
\end{lemma}
\textit{Proof}. Follows easily from the definiton. $\blacksquare$

The following proposition is an immediate generalization of
\cite{convolution}, Lemma 2.1. We include the proof because of its
importance and because it relies on the particular structure of
the measure $\mu$.

\begin{proposition} \label{prop:neighborweights}
Let $\sigma,\sigma'\in\mathcal{A}^k$, and suppose that $\sigma$
and $\sigma'$ are such that $|\pi(\sigma)-\pi(\sigma')|=d^{-k}$;
i.e. their projections are ``neighbors''. Then
\[
\frac{1}{k\theta} \le \frac{\eta(\sigma)}{\eta(\sigma')} \le
k\theta,
\]
where
\[
\theta = \frac{\max_{0\le i\le m} p_i }{\min_{0\le i\le m} p_i}.
\]
\end{proposition}
\textit{Proof}. We proceed by induction; the result is clear for
$k=1$, so assume it is valid for some $k$, and let
$\sigma,\sigma'\in\mathcal{A}^{k+1}$. Without loss of generality
suppose $\pi(\sigma') = \pi(\sigma) + d^{-(k+1)}$. Let
$A=[\pi(\sigma)]_{k+1}$, $A'=[\pi(\sigma')]_{k+1}$. In addition
let $B$ be the subset of $A$ of words ending in $m$, and note
that $\omega\mapsto (\omega_1,\ldots,\omega_{k},\omega_{k+1}+1)$
is an injective map from $A\backslash B$ into $A'$. Hence
\begin{equation} \label{eq:neighborweights:1}
\sum_{\omega\in A\backslash B} \mathbf{p}(\omega) \le
\sum_{\omega\in A\backslash B} \mathbf{p}(\omega|k) \theta
p_{\omega_{k+1}+1} \le \theta \sum_{\omega'\in A'}
\mathbf{p}(\omega') = \theta\eta(\sigma') .
\end{equation}
Now let $\omega\in B$, and write
\[
\pi(\omega) = \pi(\omega|k) + md^{-(k+1)} = jd^{-k} + md^{-(k+1)}.
\]
Since $\pi(\omega)$ is not maximal in $\pi(\mathcal{A}^{k+1})$,
$jd^{-k}$ is not maximal in $\pi(\mathcal{A}^k)$. Thus there is
$\omega'\in\mathcal{A}^k$ such that $\pi(\omega')=(j+1)d^{-k}$.
By the inductive hypothesis, $\eta(\omega|k)\le
k\theta\eta(\omega')$, and therefore
\begin{equation} \label{eq:neighborweights:2}
\sum_{\omega\in B} \mathbf{p}(\omega) = \eta(\omega|k) p_m \le
k\theta\eta(\omega')p_m \le k\theta \eta(\omega')p_{m-d+1} \le
k\theta\eta(\sigma'),
\end{equation}
since
\[
\pi(\omega',m-d+1)=(j+1)d^{-k}+(m-d+1)d^{-(k+1)}= \pi(\sigma').
\]
(it is here that we use regularity; more precisely, that $p_m\le p_{m-d+1}$). Now
combining (\ref{eq:neighborweights:1}) and
(\ref{eq:neighborweights:2}) we get
\[
\eta(\sigma) = \sum_{\omega\in A\backslash B} \mathbf{p}(\omega)
+ \sum_{\omega\in B} \textbf{p}(\omega) \le \theta \eta(\sigma')
+ k\theta\eta(\sigma') = (k+1)\eta(\sigma').
\]
The other inequality follows in the same way. $\blacksquare$

\begin{proposition} \label{prop:dimensionformula}
For all $\omega\in\mathcal{A}^{\mathbb{N}}$ we have
\[
\udim\mu(\pi(\omega)) = \limsup_{k\rightarrow\infty}
\frac{\delta\log\eta(\omega|k)}{k},
\]
and analogously for the lower dimension.
\end{proposition}
\textit{Proof}. It follows exactly like in \cite{convolution},
Proposition 2.2, using Proposition \ref{prop:neighborweights}
instead of Lemma 2.1 of \cite{convolution}. $\blacksquare$

The following result provides a formula for $\udim\mu$ (and
$\ldim\mu$) which will play a symmetric role to that of
Proposition \ref{prop:dimensionformula}.
\begin{proposition} \label{prop:dimensionmatrix}
For every $\omega\in\mathcal{A}^{\mathbb{N}}$, $\omega\ne
(0,0,\ldots)$ or $(m,m,\ldots)$, we have
\[
\udim\mu(\pi(\omega)) = \limsup_{k\rightarrow\infty}
\frac{\delta\log \|M(\omega|k)\|}{k},
\]
and analogously for $\ldim\mu(\pi(\omega))$, where $\|\cdot\|$
denotes any consistent operator.
\end{proposition}
\textit{Proof}. Fix $\omega\in\mathcal{A}^{\mathbb{N}}$,
$\omega\ne (0,0,\ldots),(m,m,\ldots)$. Define a sequence
$\{v_j:\mathbb{Z}\rightarrow\mathbb{R}\}_{j\in\mathbb{N}}$ by
$v_0(l)=\delta_{0l}=1$ if $l=0$ and $0$ otherwise; and, for $j\ge
1$,
\[
v_j(l) = \sum\{ \mathbf{p}(\sigma):\sigma\in
[\pi(\omega|j)+ld^{-j}]_j \}.
\]

Write $L=\omega_{j+1}+k$, and observe that for every
$n\in\mathbb{Z}$
\[
(\sigma,n)\in [\pi(\omega|(j+1))+kd^{-(j+1)}]_{j+1}
\Longleftrightarrow  \sigma \in [\pi(\omega|j)+(L-n)d^{-(j+1)}]_j;
\]
and this happens only if $L-n$ is a multiple of $d$ (otherwise
the sets involved are empty). Hence, recalling that $p_n=0$ if
$n\notin \{0,\ldots,m\}$,
\[
\begin{array}{lll}
v_{j+1}(k) & = &  \sum\{ \mathbf{p}(\sigma,n):(\sigma,n)\in [\pi(\omega|(j+1))+kd^{-(j+1)}]_{j+1}\} \medskip\\
& = & \sum_{l\in\mathbb{Z}} p_{(L-ld)}
\sum\{\mathbf{p}(\sigma):\sigma\in[\pi(\omega|j)+ld^{-j}]_j\} \medskip\\
& = & \sum_{l\in\mathbb{Z}} p_{(L-ld)} \,v_j(l).
\end{array}
\]
In the last sum it is enough for $l$ to run from $-a$ to $a$ whenever
$-a\le k\le a$. Indeed,
\[
\begin{array}{lll}
p_{(L-ld)} \ne 0 & \Rightarrow & 0\le
\omega_{j+1}+k-ld \le m \medskip\\
& \Rightarrow  & -m \le k-ld \le m  \medskip\\
& \Rightarrow &|l|d \le m+|k| \le m+a \medskip \\
& \Rightarrow  & |l| \le a.
\end{array}
\]
Thus we obtain
\begin{equation} \label{eq:dimensionmatrix:1}
v_{j+1}(k) = \sum_{l=-a}^a M_{\omega_{j+1}}(k,l)v_j(l) \quad
(-a\le k,l\le a).
\end{equation}
Now it follows by induction that $v_{j} = M(\omega|j) v_{0}$,
where we now consider the $v_j$ as vectors of length $2a+1$
indexed by $(-a,\ldots,a)$ rather than functions from
$\mathbb{Z}$ to $\mathbb{R}$. Recalling the definition of $v_0$
one sees that the central column of $M(\omega|j)$ is precisely
$v_j$; in particular, $M(\omega|j)(0,0)=\eta(\omega|j)$.

Now observe that if $0<\pi(\omega|j)+ld^{-j}<\xi$ then
$[\pi(\omega|j)+ld^{-j}]_j$ is nonempty, whence $v_j(l)>0$. Since
$\omega\neq (0,0,\ldots),(m,m,\ldots)$, this shows that there is
$j_0$ such that $v_{j_0}$ contains no zero coordinate. Let
\[
C = \frac{\max_{-a\le k,l\le a} M(\omega|j_0)(k,l)}{\min_{-a\le
l\le a}v_{j_0}(l)}.
\]
All columns of $M(\omega|j_0)$ are bounded by $C v_{j_0}$, and
therefore it follows from (\ref{eq:dimensionmatrix:1}) that this
happens for every $j\ge j_0$. Moreover, Proposition
\ref{prop:neighborweights} implies that
\[
\sum_{l=-a}^a v_j(l) \le (2a+1)(\theta j)^a \eta(\omega|j).
\]
Now let $\|\cdot\|_1$ denote the 1-norm on $\mathbb{R}^{2a+1}$.
Since all norms are equivalent, it suffices to prove the
proposition for the associated operator norm, which is known to be
the maximum of the 1-norms of the columns. We have
\[
\eta(\omega|j) \le  \|v_j\|_1 \le \|M(\omega|j)\|_1 \le  C(2a+1)
 (\theta j)^a\eta(\omega|j).
\]
Now it is enough to take logarithms, divide by$j$ and recall
Proposition \ref{prop:dimensionformula} to complete the proof.
$\blacksquare$

We record an interesting fact that emerged in the previous proof.

\begin{corollary} \label{coro:estimatenormeta}
Let $K=[a,b]$, with $0<a<b<\xi$. There are constants $C=C(K), D$
such that if $\sigma\in\ma^*$ and $\pi(\sigma)\in K$, then
\[
\|M(\sigma)\|_1 \le C |\sigma|^{D} \eta(\sigma),
\]
where $\|\cdot\|_1$ denotes the $1$-operator norm.
\end{corollary}
\textit{Proof}. The number of consecutive digits $0$ or $m$ at
the beginning of $\sigma$ is clearly bounded by a constant
depending only on $K$. Hence the result follows easily from the
proof of the proposition. $\blacksquare$

We indicate that while $\eta$ is supermultiplicative (in the sense
$\eta(\sigma,\sigma')\ge\eta(\sigma)\eta(\sigma')$), consistent norms
are submultiplicative; this fact will be strategically used
throughout the paper. We underline that the preceding lemma is
false for $\omega=(0,0,\ldots)$ if $p_0$ is strictly less than the
intermediate weights. The reason is that in this case $v_0$ has no
component in the direction of any eigenvector of $M_0$
corresponding to the Perron eigenvalue, so $\|M_0^k\|$ cannot be
compared with $\|M_0^k v_0\|$. The same is true of
$\omega=(m,m,\ldots)$ if $p_m<p_i$ for $i=1,\ldots,m-1$. This is
another way to look at the distinguished role that the points $0$
and $\xi$ play in the multifractal analysis of $\mu$.

Periodic sequences provide the simplest example of points where
the local dimension exists and can be computed, as the next lemma
shows.

\begin{lemma} \label{lemma:periodic}
Let $\omega\in\mathcal{A}^{\mathbb{N}}$ be periodic with period
$\sigma\in\mathcal{A}^k$ ($\sigma\ne (0,\ldots,0)$ or
$(m,\ldots,m)$). Then the local dimension of $\mu$ at
$\pi(\omega)$ exists and is given by
\[
\dim\mu(\pi(\omega)) = \delta k^{-1}\log\rho(M(\sigma)).
\]
\end{lemma}
\textit{Proof}. Let $\|\cdot\|$ be a consistent norm. Note that
$M(\omega|kj)=M(\sigma)^j$; hence (\ref{eq:radiusnorm}) shows that
$\|M(\omega|kj)\|^{1/j} \rightarrow \rho(M(\sigma))$, and
therefore
\begin{equation} \label{eq:periodic:1}
\lim_{j\rightarrow\infty} \frac{\delta \log \|M(\omega|kj)\|}
{kj} = \frac{\delta \log\rho(M(\sigma))}{k}.
\end{equation}
Now let $q=kj+r$, $0\le r<k$. Using the consistency of
$\|\cdot\|$ and setting $C=\max\{
\|M(\tilde{\sigma})\|:|\tilde{\sigma}|<k\}$ we get
\[
(1/C) \|M(\omega|k(j+1))\|\le \|M(\omega|q)\|\le C
\|M(\omega|kj)\|,
\]
This, together with (\ref{eq:periodic:1}), shows that
\[
\lim_{j\rightarrow\infty} \frac{\delta \log\|M(\omega|j)\|}{j} =
\frac{\delta \rho(M(\sigma))}{k},
\]
which completes the proof. $\blacksquare$

\begin{corollary} \label{coro:inequalityetarho}
$\eta(\sigma) \le \rho(M(\sigma))$ for every $\sigma\in\ma^*$.
\end{corollary}
\textit{Proof}. Let $\omega\in\ma^\bbn$ be periodic with period
$\sigma$. Since, by Lemma \ref{lemma:supermultiplicativity},
$j\log\eta(\sigma) \le \log\eta(\omega|j\sigma)$, we obtain
\[
\frac{\delta \log\rho(M(\sigma))}{|\sigma|} =
\lim_{j\rightarrow\infty} \frac{\delta
\log\eta(\omega|j|\sigma|)}{j|\sigma|} \ge \frac{\delta\log\eta(\sigma)}{|\sigma|}.
\]
From here the corollary follows immediately. $\blacksquare$

The following lemma illustrates the main advantage of the
existence of barrier digits; compare with Lemma
\ref{lemma:supermultiplicativity}.

\begin{lemma} \label{lemma:barrierdigit}
For all $\sigma,\sigma'\in\mathcal{A}^*$ for which $\sigma'_1$ is
a barrier digit, we have
\[
\eta(\sigma,\sigma') = \eta(\sigma)\eta(\sigma').
\]
\end{lemma}
\textit{Proof}. Take $\omega,\omega'\in\mathcal{A}^*$ such that
$|\omega|=|\sigma|,|\omega'|=|\sigma'|$ and $\pi(\omega,\omega')=
\pi(\sigma,\sigma')$. Therefore
\[
|\pi(\omega')-\pi(\sigma')|= d^{|\omega|} |\pi(\omega)-\pi(\sigma)|.
\]
Hence either $|\pi(\omega')-\pi(\sigma')|\ge 1$ or $\pi(\omega')=\pi(\sigma')$. The first is impossible since $\sigma'_1$ is a barrier digit, so we must have $\pi(\omega')=\pi(\sigma')$ and $\pi(\omega)=\pi(\sigma)$. But then
\begin{eqnarray}
\eta(\sigma,\sigma') & = & \sum \{ \mathbf{p}(\omega)\mathbf{p}(\omega') : \pi(\omega,\omega')=\pi(\sigma,\sigma'),
|\omega|=|\sigma|,|\omega'|=|\sigma'|\} \nonumber\\
& = & \sum\{ \mathbf{p}(\omega):\pi(\omega)=\pi(\sigma)  \}\sum\{ \mathbf{p}(\omega'):\pi(\omega')=\pi(\sigma')  \} \nonumber\\
& = & \eta(\sigma)\eta(\sigma'). \nonumber\,\,\blacksquare
\end{eqnarray}

Recall (\ref{eq:shat}) and note
that $\hat{S}_k$ is a strictly decreasing function such that
$\hat{S}_k(0)=\#\Xi_k>1$ and $\hat{S}_k(1) <
\sum_{\sigma\in\Xi_k} \eta(\sigma) = 1$. Hence there is a unique
``auxiliary exponent'' $0<\beta_k<1$ such that
$\hat{S}_k(\beta_k)=1$.
\begin{lemma} \label{lemma:auxiliaryexponent}
$\lim_{k\rightarrow\infty} \beta_k = 1$.
\end{lemma}
\textit{Proof}. Using Lemma \ref{lemma:supermultiplicativity}
once more we get $\eta(b,\sigma)\ge p_b\eta(\sigma)$. Hence
\begin{equation} \label{eq:auxexp1}
1 \ge p_b^{\beta_k} \sum_{\sigma\in\Xi_{k-1}}
\eta(\sigma)^{\beta_k} > p_b \max_{\sigma\in\Xi_{k-1}}
\eta(\sigma)^{\beta_k-1},
\end{equation}
where for the second inequality we used that $\sum_{\sigma\in\Xi_{k-1}}\eta(\sigma)=1$. Note that
\begin{equation} \label{eq:auxexp2}
\lim_{k\rightarrow \infty} \,\max_{\sigma\in\Xi_{k-1}} \eta(\sigma) = 0.
\end{equation}
There are several ways to see this. For instance, we know from the proof of Proposition \ref{prop:dimensionmatrix} that
\[
\eta(\sigma) \le \|M(\sigma)\|_1 \le (\max \{\|M(i)\|_1\})^{|\sigma|},
\]
where $\|M(i)\|_1$ is the maximum of the $1$-norms of the columns of $M_i$, whose non-zero coordinates are some, but not all, of the $p_i$. Hence $\eta(\sigma) < c^{|\sigma|}$ for some $c<1$, and this establishes (\ref{eq:auxexp2}) which, together with (\ref{eq:auxexp1}), imply the lemma.
$\blacksquare$

\section{Proof of the main results}

\textbf{Proof of Theorem \ref{theo:extremedimensions}}. Let
$\{\omega^j\}$ be a sequence in $\mathcal{A}^{\mathbb{N}}$ such
that $\ldim\mu(\pi(\omega^j)) \rightarrow
\inf\underline{\Delta}\mu$ as $j\rightarrow\infty$. For each $j$
choose $k_j$ such that
\begin{equation} \label{eq:extremedimensions1}
\frac{\delta\log \eta(\omega^j|k_j)}{k_j} <
\ldim\mu(\pi(\omega^j))+1/j.
\end{equation}
Let $\omega'^j$ be the periodic sequence with period
$\omega^j|k_j$. It follows from Lemma \ref{lemma:periodic} and Corollary \ref{coro:inequalityetarho} that $\dim_\mu(\pi(\omega'_j))$ exists and
\[
\dim_\mu(\pi(\omega'_j)) = \frac{\delta\log \rho(M(\omega^j|k_j))}{k_j} \le \frac{\delta\log \eta(\omega^j|k_j)}{k_j}.
\]
Together with (\ref{eq:extremedimensions1})this proves that $\inf \Delta\mu =\inf
\underline{\Delta}\mu$.

Now let $k \in \mathbb{N}$ and $M= M(\sigma)$ for some
$\sigma\in\mathcal{A}^k$. Let $\omega\in\mathcal{A}^{\mathbb{N}}$
be periodic with period $\sigma$. It follows from Lemma
\ref{lemma:periodic} that $\inf\Delta\mu\le \delta
k^{-1}\rho(M)$. Hence
\[
\inf\Delta\mu\le
\delta\log\tilde{\rho}(M_0,\ldots,M_m).
\]
For the other inequality fix $\varepsilon>0$. The previously stated facts about
periodic sequences and Lemma \ref{lemma:periodic} show the
existence of $\sigma\in\mathcal{A}^*$ such that
\[
\delta|\sigma|^{-1}\log\rho(M(\sigma)) <
\inf\Delta\mu+\varepsilon.
\]
Since $|\sigma|$ can be arbitrarily large we conclude that
\[
\delta
\log\tilde{\rho}(M_0,\ldots,M_m) \le \inf\Delta\mu.
\]
From Proposition \ref{prop:dimensionformula} it immediately follows
 that
\[
\dim\mu(0) = \delta\log p_0; \quad  \dim\mu(\xi)= \delta\log p_m.
\]
Since $\eta(\sigma)\ge
\mathbf{p}(\sigma)$ and $\mathbf{p}$ is regular, $\udim\mu(x)\le \overline{\alpha}$ for all
$x\in [0,\xi]$.

To show that $\overline{\alpha}$ is isolated in
$\overline{\Delta}\mu$ if $p_0<p_i$ for $i=1,\ldots,m-1$, we can
proceed exactly like in \cite{convolution}, Theorem 1.1; we
therefore omit the details.

It remains to prove (\ref{eq:estimatedeltamu}). It follows from Proposition \ref{prop:dimensionmatrix} and submultiplicativity of consistent norms that
$\udim\mu(\pi(T\sigma))\ge \udim\mu(\pi(\sigma))$. Note that
every $x\in(0,\xi)$ can be represented by a sequence containing
at least one digit other than $0$ and $m$ (it has to be different
from $(0,0,\ldots)$ and $(m,m,\ldots)$; replace any occurrence of
$(0,m)$ by $(1,m-d)$ and any occurrence of $(m,0)$ by $(m-1,d)$).
Therefore by taking such a sequence and shifting it we see that
\[
\sup (\Delta\mu\backslash\{ \dim\mu(0),\dim\mu(\xi) \}) =
\sup\{\udim\mu(\pi(\omega)):\omega\in\mathcal{A}^{\mathbb{N}},\omega_1\notin\{0,m\}\}.
\]
Now we can proceed as above, but now using Proposition
\ref{prop:dimensionmatrix} instead of Proposition
\ref{prop:dimensionformula}, to show that in the supremum above we
can restrict ourselves to periodic sequences. The result then
follows from Lemma \ref{lemma:periodic}. $\blacksquare$

The formulae for $\underline{\alpha}$ and $\alpha^*$ given in
Theorem \ref{theo:extremedimensions} are not very useful unless
there is some way to estimate their value in concrete cases. This
is a difficult problem, but we nevertheless have the following
corollary.
\begin{corollary} \label{coro:approximationalpha}
Fix a consistent matrix norm $\|\cdot\|$, and let
\[
\begin{array}{lll}
\tilde{\rho}_k & = &
\max\{\rho(M(\sigma)):\sigma\in\mathcal{A}^k\}^{1/k};\smallskip\\
\hat{\rho}_k & = &
\max\{\|M(\sigma)\|:\sigma\in\mathcal{A}^k\}^{1/k};\smallskip\\
\tilde{\rho}^*_k & = &
\min\{\rho(M(\sigma)):\sigma\in\mathcal{A}^k,\sigma_1 \notin \{0,m\}\}^{1/k};\smallskip\\
\hat{\rho}^*_k & = & \min\{\eta(\sigma):\sigma\in\mathcal{A}^k,
\sigma_1 \notin \{0,m\}\}^{1/k}.
\end{array}
\]

Then
\[
\begin{array}{lllll}
\delta\log\hat{\rho}_k & \le &  \underline{\alpha} & \le & \delta
\log\tilde{\rho}_k;\medskip\\
\delta\log\tilde{\rho}^*_k & \le &  \alpha^* & \le & \delta
\log\hat{\rho}^*_k.
\end{array}
\]
for every $k$. Moreover,
\[
\begin{array}{lllll}
\lim_{k\rightarrow\infty} \delta\log\hat{\rho}_k & = &
\lim_{k\rightarrow\infty} \delta\log\tilde{\rho}_k & = &
\underline{\alpha};\medskip\\
\lim_{k\rightarrow\infty} \delta\log\hat{\rho}^*_k & = &
\lim_{k\rightarrow\infty} \delta\log\tilde{\rho}^*_k & = &
\alpha^*.
\end{array}
\]
\end{corollary}
\textit{Proof}. Let $\mathcal{M} =\{M_0,\ldots,M_m\}$. The
inequalities $\tilde{\rho}_k\le \tilde{\rho}(\mathcal{M})\le
\hat{\rho}_k$ follow easily from (\ref{eq:radiusnorm}). The fact
that $\hat{\rho}_k-\tilde{\rho}_k\rightarrow 0$ as
$k\rightarrow\infty$ was proved in \cite{boundedmatrices} for any
bounded set of matrices; hence
\[
\lim_{k\rightarrow\infty} \delta\log\hat{\rho}_k  =
\lim_{k\rightarrow\infty} \delta\log\tilde{\rho}_k  =
\underline{\alpha}.
\]

Now we turn to the approximations of $\alpha^*$. From
(\ref{eq:radiusnorm}), Corollary \ref{coro:estimatenormeta} and
Corollary \ref{coro:inequalityetarho} we obtain
\[
\eta(\sigma) \le \rho(\sigma) \le C|\sigma|^D \eta(\sigma),
\]
for every $\sigma\in\ma^*$ such that $\sigma_1\notin \{0,m\}$,
where $C$ and $D$ are independent of $\sigma$. From here it
follows that
\[
\delta\log\tilde{\rho}^*_k  \le  \alpha^*  \le  \delta
\log\hat{\rho}^*_k \le \delta \log\tilde{\rho}^*_k - \delta
\frac{\log(C k^D)}{k}.
\]
In particular,
\[
\lim_{k\rightarrow\infty} \delta\log\hat{\rho}^*_k  =
\lim_{k\rightarrow\infty} \delta\log\tilde{\rho}^*_k  = \alpha^*. \quad\blacksquare
\]

\bigskip

\textbf{Proof of Proposition \ref{prop:lqspectrum}}. Let $\mu_k$
be the discrete measure assigning mass $\eta(\sigma)$ to
$\pi(\sigma)$ for $\sigma\in\mathcal{A}^k$. Note that, for $x\in\supp\mu$,
\begin{equation} \label{eq:prop2gen:1}
\mu(B(x,d^{-k})) \le \mu_k(B(x,cd^{-k})) \le \mu(B(x,2cd^{-k})),
\end{equation}
where $c=1+\diam(\supp\mu)=1+\xi$. Let $x=\pi(\omega)$, and
observe that $\pi(\omega|k)\in B(x,cd^{-k})$. On the other hand,
$B(x,cd^{-k})$ contains no more than $\lfloor 2c \rfloor$
consecutive points of $\pi(\mathcal{A}^k)$. Hence, by Proposition
\ref{prop:neighborweights},
\begin{equation} \label{eq:prop2gen:2}
\eta(\omega|k) \le \mu_k(B(x,cd^{-k})) \le \lfloor 2c \rfloor
(\theta k)^{\lfloor 2c \rfloor} \eta(\omega|k).
\end{equation}
From (\ref{eq:prop2gen:1}) and (\ref{eq:prop2gen:2}) a routine,
but maybe a little bit tedious, calculation shows that
\begin{equation} \label{eq:expressionlqspectrum}
\underline{\tau}(q) = \liminf_{k\rightarrow\infty}
\frac{\delta\log \overline{S}_k(q)}{k},
\end{equation}
We will now show that the limit in (\ref{eq:expressionlqspectrum})
exists. Write
\[
\Xi^*_{n,k}=\{(\sigma,\omega):\sigma\in
\Xi_n,\omega\in\Xi_k\}\subset\mathcal{A}^{k+n}.
\]
Note that $\pi(\Xi^*_{n,k})=\pi(\mathcal{A}^{k+n})$.
On the other hand, if $\pi(\sigma,\omega)=\pi(\sigma',\omega')$
for $\sigma,\sigma'\in\Xi_n; \omega,\omega'\in\Xi_k$, then
\[
|\pi(\sigma)-\pi(\sigma')| = d^{-n} |\pi(\omega)-\pi(\omega')| \le  d^{-n}\xi.
\]
Hence at most $\lfloor \xi\rfloor+1$ elements of  $\Xi^*_{n,k}$
project onto the same number. Using once again Lemma
\ref{lemma:supermultiplicativity} we obtain
\[
\overline{S}_n(q)\overline{S}_k(q)\le
\sum_{\sigma\in\Xi_n,\tau\in\Xi_k}\eta(\sigma,\tau)^q \le (\lfloor
\xi\rfloor+1) \overline{S}_{n+k}(q),
\]
for $q\ge 0$, and
\[
\overline{S}_n(q)\overline{S}_k(q)\ge \sum_{ \sigma \in \Xi_n,
\tau \in \Xi_k} \eta(\sigma,\tau)^q \ge \overline{S}_{n+k}(q),
\]
for $q\le 0$. In either case, sub/supermultiplicativity shows that the limit in
(\ref{eq:expressionlqspectrum}) exists. $\blacksquare$

Two remarks are in order. First, Peres and Solomyak
\cite{lqspectrum} have proved that for \textit{any} self-similar
measure the $L^q$-spectrum exists in the range $q\ge 0$,
regardless of separation; their proof also
relies on submultiplicativity. Second, Lau and Ngai showed in
\cite{wsc} that the Legendre transform of
$\underline{\tau}_\mu(q)$ is always an upper bound for
$f_{\mu,H}(\alpha)$.

\bigskip

We begin now our investigation of the case $m<2d-2$. Recall that
by iterating the IFS if necessary we can assume without loss of generality the
existence of a barrier digit $b$.

\bigskip

\textbf{Proof of theorem \ref{theo:deltamupart}}. It is clear that
$\dim \mu(0)=\delta\log p_0$ and $\dim \mu(\xi)=\delta \log p_m$.
Hence it suffices to prove that $(\underline{\alpha},\alpha^*)
\subset \Delta\mu$.

Fix $k\in\mathbb{N}$, and let $\mu^k$ be the attractor of the IFS
\begin{equation} \label{eq:ifsbarrierdigit}
\left\{\left(d^{-k}(x+d^k\pi(b,\sigma)),\eta(b,\sigma)^{\beta_k}\right):
\sigma\in\Xi_{k-1}\right\},
\end{equation}
where $\beta_k$ is the auxiliary exponent defined before Lemma \ref{lemma:auxiliaryexponent}.

Note that, since $b$ is a barrier, $0<\pi(b,\sigma)<1$. Thus the
maps are of the form $d^{-k}(x+j)$, with $0<j<d^k$, and the IFS
verifies the strong separation condition. This allows us to use
the multifractal theory developed in \cite{mauldin}; see also
\cite{falconer3}, Chapter 11. Let
\[
\mathcal{A}_k= \{(b,\sigma):\sigma\in\Xi_{k-1}\}.
\]
Endow $\mathcal{A}^{\mathbb{N}}_k$ with the Bernoulli measure $\nu^k$ for
the weights $\eta(b,\sigma)^{\beta_k}$, and denote by
$\pi^k:\mathcal{A}^{\mathbb{N}}_k\rightarrow \mathbb{R}$ the
canonical projection. Observe that
\begin{equation} \label{eq:udimstrongseparation}
\udim\mu^k(\pi^k(\omega)) = \limsup_{j\rightarrow\infty}
\frac{\delta\log\nu^k([\omega|j])}{jk},
\end{equation}
and analogously for the lower dimension, where $[\omega|j]$
denotes the cylinder generated by $\omega|j$. There is a
canonical map
$\phi^k:\mathcal{A}^{\mathbb{N}}_k\rightarrow\mathcal{A}^{\mathbb{N}}$;
clearly $\pi\circ\phi^k=\pi^k$. Moreover,
\[
\nu^k([\omega|j])=\eta(\omega_1)^{\beta_k} \cdots \eta(\omega|j)^
{\beta_k}= \eta(\phi^k(\omega)|kj)^{\beta_k},
\]
where for the last equality we used Lemma
\ref{lemma:barrierdigit}. From this,
(\ref{eq:udimstrongseparation}) and Proposition
\ref{prop:dimensionformula} it follows that if $x=\pi^k(\omega)$
for some $\omega\in\mathcal{A}^{\mathbb{N}}_k$, then
\begin{equation} \label{eq:dimensionequality}
\udim\mu^k(x)=\beta_k \udim\mu(x),
\end{equation}
and analogously for the lower dimension. This equality is the key to estimating the multifractal spectrum of $\mu$, by reducing it to the study of the spectrum of the measures $\mu^k$ (which are well-understood). As a first instance of this, note that
$\Delta\mu\supset \beta_k^{-1}\Delta{\mu^k}$. Hence we obtain
(see \cite{mauldin}) $\Delta\mu \supset [ \underline{\alpha}_k,
\overline{\alpha}_k]$, where
\[
\underline{\alpha}_k=\frac{\delta\max_{\sigma\in\Xi_{k-1}}\log\eta(b,\sigma)}{
\beta_k k};
\]
\[
\overline{\alpha}_k=\frac{\delta\min_{\sigma\in\Xi_{k-1}}\log\eta(b,\sigma)}{
\beta_k k}.
\]
Letting $k$ run through the positive integers we get $\Delta\mu
\supset (\inf_k \underline{\alpha}_k,
\sup_k\overline{\alpha}_k)$. To complete the proof we will now
show that $\inf_{k} \underline{\alpha_k} = \underline{\alpha}$ and
$\sup_k \overline{\alpha}_k = \alpha^*$. We have that
$\underline{\alpha}\le \underline{\alpha}_k$ for all $k$; hence
it suffices to show that $\liminf_{k\rightarrow\infty}
\underline{\alpha}_k\le \underline{\alpha}$ . Since
$\eta(b,\sigma)\ge p_b \eta(\sigma)$ we have
\[
\underline{\alpha}_k \le \frac{\delta\log
p_b+\delta\max_{\sigma\in\Xi_{k-1}}\log \eta(\sigma)}{\beta_k k}.
\]
Thus, using Lemma \ref{lemma:auxiliaryexponent},
\[
\liminf_{k\rightarrow\infty}\underline{\alpha}_k =
\liminf_{k\rightarrow\infty}\delta k^{-1}\max_{\sigma\in \Xi_{k}}
\log \eta(\sigma) \le \underline{\alpha}.
\]
For the other equality, $\sup_k \overline{\alpha}_k=\alpha^*$,
observe that, since
$\udim\mu(\pi(b,\sigma))\ge\udim\mu(\pi(\sigma))$ (by Proposition
\ref{prop:dimensionmatrix}),
\[
\alpha^* = \sup_\omega \udim\mu(\pi(b,\omega)) \le \delta \sup_k
k^{-1} \min_{\sigma\in\Xi_{k-1}} \log\eta(b,\sigma) \le \sup_k
\overline{\alpha}_k.
\]
This concludes the proof. $\blacksquare$

\bigskip

\textbf{Proof of Theorem \ref{theo:lqspectrumpart}}. To make the
proof easier to read we will split it into three parts.

\bigskip

\textbf{First step}. We begin by showing that the limit in the
definition of $\hat{\tau}$ exists. Note that if we write
\[
\Xi_{k,n}^* = \{(\sigma,b,\sigma'):\sigma\in\Xi_{k-1},\sigma'\in
\Xi_{n-1}\},
\]
then $\pi$ restricted to $\Xi_{k,n}^*$ is injective (this follows from
$b$ being a barrier digit), and
\begin{equation} \label{eq:comparisonxi}
\{ \eta(b,\omega):\omega\in \Xi_{k,n}^* \} \subset \{ \eta(b,\sigma):\sigma\in \Xi_{k+n-1} \}.
\end{equation}
Hence we have
\[
\begin{array}{lll}
\hat{S}_k(q)\hat{S}_n(q) & = & \sum\{
\eta(b,\sigma,b,\sigma')^q: \sigma\in\Xi_k,\sigma'\in\Xi_n \} \medskip \\
& = & \sum\{\eta(b,\omega)^q:\omega\in\Xi^*_{k,n}\} \medskip \\
& \le & \hat{S}_{k+n}(q),
\end{array}
\]
where for the first equality we again used Lemma
\ref{lemma:barrierdigit}, and for the last inequality we used (\ref{eq:comparisonxi}).

 Therefore the sequence
$\log\hat{S}_{n}(q)$ is superadditive, whence the limit of
$a_k=\delta k^{-1} \log \hat{S}_k(q)$ exists and is equal to
$\sup_k a_k$.

\bigskip

\textbf{Second step}. The estimates used in establishing
Proposition \ref{prop:lqspectrum} also show the following: if $K$
is any closed subinterval of $[0,\xi]$ and
\[
S^K_j(q) = \sum\{\eta(\sigma)^q:\sigma\in\Xi_k \textrm{ and
}\pi(\sigma) \in K\},
\]
then
\[
\underline{\tau}^K(q) = \liminf_{k\rightarrow\infty} \delta k^{-1}
\log S^K_j(q),
\]
where $\underline{\tau}^K$ denotes the lower $L^q$-spectrum of
$\mu|_K$; and an analogous assertion holds for the upper limit.

Choose any $\sigma\in\ma^*$ such that for any $\omega\in\ma^\bbn$
beginning with $\sigma$, $\pi(\omega)\in K$. Let $k=|\sigma|$.
For all $j>k+1$ we have
\[
\begin{array}{lll}
S^K_{j}(q) & \ge &  \sum\{\eta(\sigma,b,\sigma')^q:\sigma'\in
\Xi_{j-k-1}\} \medskip \\
& = & \eta(\sigma)^q \hat{S}_{j-k}(q),
\end{array}
\]
by Lemma \ref{lemma:barrierdigit}. From here we obtain
\begin{equation} \label{eq:lqspectrumpart:step2:0}
\limsup_{k\rightarrow\infty} \delta k^{-1}\log S^K_{j}(q) \le
\hat{\tau}(q).
\end{equation}
For the opposite inequality we consider the cases $q\ge 0$ and
$q<0$ separately. For nonnegative $q$ we have $\eta(\sigma)^q \le
p_b^{-q}\eta(b,\sigma)^q$, whence
\begin{equation} \label{eq:lqspectrumpart:step2:1}
S^K_{j}(q) \le p_b^{-q} \hat{S}_j(q).
\end{equation}
For negative $q$ we use Corollary \ref{coro:estimatenormeta}; at
this point we need to assume that $K$ is bounded away from $0$ and
$\xi$. From the Corollary we get
\[
\begin{array}{lll}
\eta(b,\sigma) & \le & \|M(b,\sigma)\| \medskip\\
& \le & \|M_b\| \|M(\sigma)\| \medskip\\
& \le & \|M_b\| C |\sigma|^D \eta(\sigma),
\end{array}
\]
where $C$ and $D$ do not depend on $\sigma$ (nor $|\sigma|$).
Therefore for $q<0$ we obtain
\[
\eta(\sigma)^q \le C' |\sigma|^{-qD} \eta(b,\sigma)^q,
\]
where $C'$ does not depend on $\sigma$. Hence
\begin{equation} \label{eq:lqspectrumpart:step2:2}
S^K_j(q) \le C' j^{-qD} \hat{S}_j(q).
\end{equation}
From (\ref{eq:lqspectrumpart:step2:1}) and
(\ref{eq:lqspectrumpart:step2:2}) we conclude
\[
\liminf_{j\rightarrow\infty} \delta j^{-1}\log S^K_{j}(q) \ge
\hat{\tau}(q).
\]
This together with (\ref{eq:lqspectrumpart:step2:0}) shows that
the $L^q$ spectrum of $\mu|_K$ exists and is equal to
$\hat{\tau}$, as desired.

Before proceeding to the next step, let us remark that the fact
that $\hat{\tau}$ is the $L^q$-spectrum of $\mu|_K$ guarantees
that $\hat{\tau}$ is a concave (hence continuous) increasing
function, and that $\hat{\tau}(0)=-1$; see \cite{wsc},
Proposition 2.3 (all these facts can also be checked directly from
the definition).

\bigskip

\textbf{Third step}. We will now prove
(\ref{eq:formulalqspectrum}). The case $q\ge 0$ follows like in
the previous step, so we will assume $q<0$. Let
\[
\begin{array}{llllll}
K_1 & = & [0,d^{-1}]; & K'_1 & = & [\xi-d^{-1},\xi];\\
K_2 & = & [d^{-1},1]; & K'_2 & = & [\xi-1,\xi-d^{-1}];\\
K_3 & = & [d^{-1},\xi-d^{-1}]; & K_4 & = & K_2 \cap K'_2.
\end{array}
\]
We will keep using the notation $S^K_j(q)$ (even if $K$ contains
$0$ or $\xi$). Let us adopt the convention that $S^K_0(q)=1$. Note
that
\begin{equation} \label{eq:lqspectrumpart:1}
\overline{S}_j(q) = S^{K_1}_j(q) + S^{K_3}_j(q) + S^{K'_1}_j(q)
\end{equation}
(there is a minor issue with the points $1$ and $\xi-1$, but this
does not affect the argument). Moreover, if $d^{-(k+1)}\le
\pi(\sigma) < d^{-k}$ for $\sigma\in\ma^*$, then $\sigma$ begins
with at least $k$ zeros, whence $\eta(\sigma) = p_0^k
\eta(\sigma')$ for $\sigma'=T^k\sigma$. Moreover,
$|\sigma'|=\sigma-k$ and $d^{-1}\le\pi(\sigma')< 1$. Hence
\[
S^{K_1}_j(q) = \sum_{k=1}^j p_0^{kq} S^{K_2}_{k-j}(q).
\]
Analogously,
\[
S^{K'_1}_j(q) = \sum_{k=1}^j p_m^{kq} S^{K'_2}_{k-j}(q).
\]
Noting that $p_m^{kq} \le p_0^{kq}$, $K_2\cup K'_2 = K_3$ and
$S_j^K(q)$ is increasing in $K$, we obtain from
(\ref{eq:lqspectrumpart:1}) that
\begin{equation} \label{eq:lqspectrumpart:2}
\sum_{k=0}^j p_0^{kq} S^{K_4}_{j-k}(q) \le \overline{S}_j(q) \le 2
\sum_{k=0}^j p_0^{kq} S^{K_3}_{j-k}(q).
\end{equation}
Fix now $\varepsilon>0$. Recall from the previous step that the $L^q$-spectrum of $\mu|_{K_i} (i=3,4)$ exists and it is given by $\hat{\tau}$; moreover, in the course of the proof we showed that
\[
\hat{\tau}(q) = \lim_{k\rightarrow\infty} \delta k^{-1} \log S_k^{K_i}(q).
\]
Recalling that $\delta = -1/\log d$ we obtain, after taking exponentials, that
\[
d^{-\hat{\tau}(q)} = \exp(-\hat{\tau}(q) \log d) = \lim_{k\rightarrow\infty} \exp\left(\frac{\log S_k^{K_i}(q)}{k}\right) = \lim_{k\rightarrow\infty}  \left(S_k^{K_i}(q)\right)^{1/k} .
\]
Therefore we see that there is a constant $C=C(\varepsilon)>0$ such that
\[
C^{-1} d^{-k(\hat{\tau}(q)+\varepsilon)} \le S^{K_4}_k(q) \le
S^{K_3}_k(q) \le C d^{-k(\hat{\tau}(q)-\varepsilon)},
\]
for all $k\ge 0$ (the second inequality follows from the fact that $K_4\subset K_3$). Plugging this into (\ref{eq:lqspectrumpart:2})
and adding we obtain
\[
C^{-1} \frac{p_0^{q(j+1)}-d^{-(j+1)(\hat{\tau}(q)+
\varepsilon)}}{p_0^{q}-d^{-(\hat{\tau}(q)+ \varepsilon)}} \le
\overline{S}_j(q) \le 2 C \frac{p_0^{q(j+1)}-d^{-(j+1)(\hat{\tau}(q)-
\varepsilon)}}{p_0^{q}-d^{-(\hat{\tau}(q)- \varepsilon)}}.
\]
Taking logarithms we get
\[
C' \left( p_0^{q(j+1)}-d^{-(j+1)(\hat{\tau}(q)+
\varepsilon)}\right) \le \overline{S}_j(q)  \le C'' \left( p_0^{q(j+1)}-d^{-(j+1)(\hat{\tau}(q)-
\varepsilon)} \right)
\]
for some positive $C',C''$ independent of $j$.

Assume first than $p_0^q < d^{\hat{\tau}(q)}$. If $\varepsilon$ is so small that $p_0^q < d^{-(\hat{\tau}(q)+\varepsilon)}$ then, by taking logarithms and then the limit as $j\rightarrow\infty$ (while keeping $\varepsilon$ fixed) we deduce
\[
\hat{\tau}(q)- \varepsilon \le \liminf_{j\rightarrow\infty} \frac{\log \overline{S}_j(q)}{-j \log d} \le  \limsup_{j\rightarrow\infty} \frac{\log \overline{S}_j(q)}{-j \log d} \le \hat{\tau}(q)+ \varepsilon.
\]
Recalling Proposition \ref{prop:lqspectrum} and letting $\varepsilon\rightarrow 0$ we obtain $\hat{\tau}(q) = \tau(q)$ in this case.

Analogously, if $p_0^q < d^{\hat{\tau}(q)}$ then, recalling that $\overline{\alpha}=\log(p_0)/(-\log d)$ or $p_0=d^{-\overline{\alpha}}$, we get
\[
\overline{\alpha} q = \lim_{j\rightarrow\infty} \delta j^{-1}\log \overline{S}_j(q),
\]
or, in other words, $\tau(q) = \overline{\alpha} q$. In short, we have
\begin{equation} \label{eq:formulalqspectrum1}
\tau(q) = \left\{
\begin{array}{lll}
\overline{\alpha} q & \textrm{ if } & \overline{\alpha} q > \hat{\tau}(q) \\
\hat{\tau}(q) & \textrm{ if } & \overline{\alpha} q \le
\hat{\tau}(q)
\end{array}
\right..
\end{equation}

(The continuity of $\hat{\tau}$ guarantees that the formula above
is also valid when $\hat{\tau}(q)=\overline{\alpha}q$). Assume
first that $\alpha^*<\overline{\alpha}$. Since, by Theorem \ref{theo:extremedimensions},
$\log\eta(b,\sigma) < k\alpha^*$ if $\sigma\in\Xi_{k-1}$, and the
number of terms in the sum $\hat{S}_k(q)$ is bounded by $C
d^{k}$, we obtain that
\[
\hat{S}_k(q) \le C d^k \exp(qk\alpha^*)\quad\textrm (q<0),
\]
and from here it follows that $\hat{\tau}(q) \ge
q\alpha^*-1$ for negative $q$. Since $\hat{\tau}(0)=-1<0$, the
concave curve $\hat{\tau}(q)$ meets the line $\overline{\alpha}q$
at a single negative point $q_0$, so (\ref{eq:formulalqspectrum})
is verified.

It remains to handle the case $\alpha^*=\overline{\alpha}$. We
have shown in the proof of Theorem \ref{theo:deltamupart} that
\[
\alpha^* = \delta \inf\{ k^{-1}
\log\eta(b,\sigma):\sigma\in\Xi_{k-1} \},
\]
whence $\hat{\tau}(q) \ge q\alpha^*=q\overline{\alpha}$ for all
$q$. Recalling (\ref{eq:formulalqspectrum1}) we see that in this
case $\tau(q)=\hat{\tau}(q)$ for all $q\in\mathbb{R}$. The proof
is now complete. $\blacksquare$

\bigskip

\textbf{Proof of Theorem \ref{theo:multifractalformalism}}. We
continue using the notation of the proof of Theorem
\ref{theo:deltamupart}. From the multifractal theory for self
similar sets under strong separation (see \cite{falconer3},
Theorem 11.5), it follows that the $L^q$-spectrum of $\mu^k$ is
given by $T_k(q)=\hat{\tau}_k(\beta_k q)$, and the multifractal
spectrum equals the Legendre transform of $T_k(q)$. Thus
\[
f_H(\alpha) \ge \dim_H\{x:\dim\mu^k(x)=\beta_k\alpha\} =
T_k^*(\beta_k\alpha) = \hat{\tau}_k^*(\alpha).
\]
The last equality follows from the definition of Legendre
transform:
\[
T_k^*(\beta_k\alpha) = \inf_{q\in\mathbb{R}}
q\beta_k\alpha-T_k(q) = \inf_{q'\in\mathbb{R}}
q'\alpha-\hat{\tau}_k(q') = \hat{\tau}_k^*(\alpha).
\]
Recall from Theorem \ref{theo:lqspectrumpart} that $\hat{\tau}=
\inf_k \hat{\tau}_k$. Let $F=\sup \hat{\tau}_k^*$. The
subadditivity used in the first step of the proof of Theorem
\ref{theo:deltamupart} shows that $F$ can be obtained as a
monotone supremum of concave functions, and thus it is concave.
Hence
\[
\begin{array}{lll}
\hat{\tau}_k^* \le F \textrm{ for all } k & \Rightarrow &
\hat{\tau}_k \ge F^* \textrm{ for all } k \medskip\\
& \Rightarrow & \hat{\tau} \ge F^* \Rightarrow \hat{\tau}^* \le F,
\end{array}
\]
where we used that the Legendre transform is involutive (i.e.
$(g^*)^* =g)$) and order-reversing on concave functions. We
conclude that $f_H(\alpha) \ge F(\alpha) \ge
\hat{\tau}^*(\alpha)$.

Let $I = d^{-1}(b,b+\xi)$, and denote by $\hat{f}_H$ the
multifractal spectrum of $\mu|_I$. We claim that $\hat{f}_H \ge
f_H$ (the opposite inequality is obvious, but it is not what we
need). To this end, observe that if $\omega\notin
\{(0,0,\ldots),(m,m,\ldots)\}$ then
$\dim\mu(\pi(b,\omega))=\dim\mu(\pi(\omega))$ (we have previously
used both inequalities separately; they follow from Propositions
\ref{prop:dimensionformula} and \ref{prop:dimensionmatrix}). Thus
$g:(0,\xi)\rightarrow I$ defined by $g(x)=d^{-1}(x+b)$ is a
bi-Lipschitz map such that $\dim\mu(x)=\dim\mu(g(x))$. The claim
now follows from the invariance of Hausdorff dimension under
bi-Lipschitz maps.

We know from Theorem \ref{theo:lqspectrumpart} that $\hat{\tau}$
is the $L^q$-spectrum of $\mu|_{\overline{I}}$. We recall the
result of Lau and Ngai (\cite{wsc}, Theorem 4.1), that the
Legendre transform of the $L^q$-spectrum is always an upper bound
for the multifractal spectrum; we remark that they do not assume
self-similarity. Hence we are able to conclude that $f_H(\alpha)
\le \hat{\tau}^*(\alpha)$, and this completes the proof.
$\blacksquare$

\bigskip

\textbf{Proof of Proposition \ref{prop:almostsuredimension}}. Let
$g_k:\mathcal{A}^{\mathbb{N}}\rightarrow\mathbb{R}$ be given by
\[
g_k(\omega) = -\log\eta(\omega|k).
\]
Observe that
\[
\int_{{\ma}^{N}} g_k(\omega)\,d\nu = \sum_{\sigma\in\mathcal{A}^k}
-\mathbf{p}(\sigma)\log\eta(\sigma) = \sum_{\sigma\in\Xi_k}
-\eta(\sigma)\log\eta(\sigma)  \ge 0.
\]
Note that Lemma \ref{lemma:supermultiplicativity} can be restated
as $g_{n+k}(\omega)\le g_n(\omega) + g_k(T^n(\omega))$ (recall
that $T$ is the shift operator). Hence the fact that the local dimension exists and is almost everywhere constant, as well as the first equality in the proposition, follow from Kingman's subadditive ergodic
theorem \cite{subadditive} applied to the system $(\mathcal{A}^\mathbb{N},\nu_{\mathbf{p}},T, \{g_k\}_k)$.

The second equality follows analogously, by considering the functions $h_k:\mathcal{A}^{\mathbb{N}}\rightarrow\mathbb{R}$ given by
\[
h_k(\omega) = - \log\|M(\omega|k)\|,
\]
for some fixed consistent norm $\|\cdot\|$, and applying the subadditive ergodic theorem to the family $\{-h_k\}$. $\blacksquare$

\section{Examples and applications}

A consequence of Theorems \ref{theo:lqspectrumpart} and
\ref{theo:multifractalformalism} is that, in the case $m<2d-2$,
the multifractal formalism holds for $\mu$ if and only if
$\alpha^* = \overline{\alpha}$. Therefore it is of interest to
find explicit necessary and sufficient conditions for the
equality of $\alpha^*$ and $\overline{\alpha}$; the next
proposition does precisely this.

\begin{proposition} \label{prop:alphastarequalsoveralpha}
$\alpha^* = \overline{\alpha}$ if and only if $m\le 2d-2$ and
$p_0=p_i$ for some $i\in \{m-d+1,\ldots,d-1\}$.
\end{proposition}
\textit{Proof}. Assume first that $p_0=p_i$ for some $i\in
\{m-d+1,\ldots,d-1\}$ (whence, in particular, $m\le 2d-2$). In
this case $a=1$ and an inspection of the matrix $M_i$ shows that
$\rho(M_i)=p_0$ (the eigenvalues are $p_i$ and $p_0$ or $p_m$ or
both; but we are assuming that these numbers are equal). Hence
\[
\alpha^* \le \dim\mu(\pi(i,i,\ldots)) = \delta \log\rho(M_i) =
\overline{\alpha},
\]
and $\alpha^*=\overline{\alpha}$ in this case.

Assume now that $p_0 < p_i$ for $i=m-d+1,\ldots,d-1$ (if there is
any such $i$).  Without loss of generality we assume $m>d$. A
simple modification of Proposition 3.4 of \cite{convolution}
shows the following: for every $s\in (0,\xi)$ there is
$\omega\in\ma^ \bbn$ such that $\pi(\omega)=s$ and the digits $0$
and $m$ appear in $\omega$ only finitely many times. Since
shifting such $\omega$ does not change upper and lower local
dimensions, we obtain that
\[
\alpha^* \le \sup\{ \udim\mu(\pi(\omega)): \omega\in
\{1,\ldots,m-1\}^\bbn \}.
\]
We will show that $\eta(\sigma)>p_0^2$ for every $\sigma\in
\{1,\ldots,m-1\}^2$. Indeed, if $\sigma_2 \notin
\{m-d+1,\ldots,d-1\}$ let
\[
\sigma^* = \left\{
\begin{array}{lll}
(\sigma_1+1,\sigma_2-d) & \textrm{ if } & \sigma_2 \ge d
\smallskip \\
(\sigma_1-1,\sigma_2+d) & \textrm{ if } & \sigma_2 \le m-d
\end{array}
\right..
\]
Note that $\pi(\sigma)=\pi(\sigma^*)$, whence
\[
\eta(\sigma) \ge \mathbf{p}(\sigma) +  \mathbf{p}(\sigma^*) \ge 2
p_0^2.
\]
If, on the other hand, $\sigma_2 \in \{m-d+1,\ldots,d-1\}$, then
\[
\eta(\sigma) \ge \mathbf{p}(\sigma) \ge p_0
\min\{p_{m-d+1},\ldots, p_{d-1}\}.
\]
Thus in any case $\eta(\sigma)\ge \tilde{p} p_0$, where
\[
\tilde{p} = \min \{p_{m-d+1},\ldots, p_{d-1}, 2 p_0\} > p_0.
\]
Splitting $\omega|2k$ in $k$ chunks of length $2$ and using
supermultiplicativity we get $\eta(\omega|2k) \ge
\left(\tilde{p}p_0\right)^k$, whence
\[
\udim\mu(\pi(\sigma)) \le \delta\frac{\log \tilde{p}+\log p_0}{2}
< \delta\log p_0 = \overline{\alpha}.
\]
Thus $\alpha^*<\overline{\alpha}$ in this case, completing the
proof. $\blacksquare$

The case in which all the weights are equal is interesting for
several reasons. On one hand, it is the most ``purely
combinatoric'' case; on the other, the vector $\mathbf{p}$ is
extremal among all allowed weights of the same length. We
investigate this case more closely. We begin by showing that if
$m\equiv -1\mod d$, then $\mu$ is absolutely continuous.

\begin{proposition} \label{prop:abscontinuous}
Let $m=nk-1$ for some $n\ge 2$, and take $p_i=1/(m+1)$ for all
$i$. The measure $\mu$ thus obtained is absolutely continuous
with a bounded density.
\end{proposition}
\textit{Proof}. From the proof of Proposition
\ref{prop:dimensionmatrix} it follows that $a$ can be replaced by
\textit{any} integer greater or equal than $\lfloor (m-d) / (d-1)
\rfloor$. For the purposes of this proof we will take $a =
\lfloor 1+\xi \rfloor$. Recall that
\[
\mu(B(\pi(\omega),d^{-k})) \le
\mu_k(B(\pi(\omega|k),(1+\xi)d^{-k})).
\]
By the choice of $a$, the right hand side is equal to the central
column of $M(\sigma|k)$ (this also follows from the proof of
Proposition \ref{prop:dimensionmatrix}). Hence
\[
\mu(B(\pi(\omega),d^{-k})) \le \|M(\omega|k)\|_1 \le C_1
\|M(\omega|k)\|_\infty \le C_1 \left(\max_{0\le i\le m}
\|M_i\|_\infty\right)^k.
\]
for some constant $C_1$. But the $\infty$-operator norm is equal
to the maximum of the $1$-norms of the rows. Observe that any
given row of $M_i$ is of the form
$(p_{j-ad},\ldots,p_j,\ldots,p_{j+ad})$, whence at most $n$ of the
coordinates are nonzero. Since $p_j=(dn)^{-1}$ for every
$j=0,\ldots, m$, it follows that $\|M_i\|_{\infty} \le d^{-1}$ for
all $i=0,\ldots,m$. Thus we obtain that
\[
\mu(B(\pi(\omega),r)) \le C_2 \lambda(B(\pi(\omega),r))
\quad\textrm{ for every } \omega\in\ma^\bbn, r>0,
\]
for some constant $C_2$, where $\lambda$ denotes Lebesgue measure
on the line. We conclude that $\mu$ is absolutely continuous and,
moreover, $d\mu/d\lambda$ is bounded. $\blacksquare$

The next proposition deals with the case $m=d$.

\begin{proposition} \label{proposition:mequalsd}
Let $m=d$ and $p_i=1/(m+1)$ for all $i$. Then
\[
\underline{\alpha} = \frac{\log(d+1)-\log \zeta}{\log d},
\]
where $\zeta=(1+\sqrt{5})/2$ is the golden number. Also
$\alpha^*=\overline{\alpha}= \log(d+1)/\log d$.
\end{proposition}
\textit{Proof}. A calculation shows that
\[
M_0 M_1 = \frac{1}{(d+1)^2}\left(
\begin{array}{lll}
1 & 1 & 0 \\
1 & 2 & 0 \\
0 & 1 & 0
\end{array}
\right) \,\Longrightarrow\, \rho(M_0 M_1) =
\left(\frac{\zeta}{d+1}\right)^2.
\]
whence
\[
\underline{\alpha} \le \frac{\delta\log(\rho(M_0 M_1))}{2} =
\frac{\log(d+1)-\log \zeta}{\log d}.
\]
We will inductively show that for all $\sigma\in\ma^j$,
\[
\begin{array}{lll}
\sigma_j \in \{0,d\} & \Rightarrow & \eta(\sigma) \le (d+1)^{-j}
F_{j} \\
\sigma_j \in \{1,\ldots,d-1\} & \Rightarrow & \eta(\sigma) \le
(d+1)^{-j} F_{j-1},
\end{array}
\]
where $F_j$ denotes the $j$-th Fibonacci number. Indeed, this is
clear for $j=1,2$. Assume it is valid for $j=1,\ldots,n$, and let
$\sigma\in\ma^{n+1}$. Observe that if $\sigma'\in
[\pi(\sigma)]_{n+1}$ then $\sigma'_{n+1}\equiv \sigma_{n+1}\mod
d$ (just multiply $\pi(\sigma)=\pi(\sigma')$ by $d^{n+1}$). Hence
two cases arise. If $0<\sigma_{n+1}<d$ we obtain
\[
\eta(\sigma) = \eta(\sigma|n) p_{\sigma(n+1)} \le (d+1)^{-n} F_n,
\]
by the inductive hypothesis. If $\sigma_{n+1} = 0$ and $\sigma'\in
[\pi(\sigma)]_{n+1}$, it must be $\sigma'_{n+1}=0$ or $d$. In the
second case $\pi(\sigma|n)=\pi(\sigma'|n)+d^{-n}$ whence,
multiplying by $d^n$, $\sigma_n \equiv \sigma'_n +1\mod d$. It
follows that $\sigma_n$ and $\sigma_n'$ cannot be simultaneously
in $\{0,d\}$, and therefore, using the inductive hypothesis,
\[
\eta(\sigma) = p_0 \eta(\sigma|n) + p_m \eta(\sigma'|n) \le
(d+1)^{-n} (F_{n-1} + F_{n-2}) = (d+1)^{-n} F_{n}.
\]
(if no such $\sigma'$ exists it is even simpler to obtain the
needed estimate). The case $\sigma_{n+1}=d$ is handled in the same
way.

Recall that $F_j=\zeta^j-{\zeta'}^j$ for some $0<\zeta'<1$. We
conclude from Proposition \ref{prop:dimensionformula} that
\[
\ldim\mu(\pi(\omega)) \ge  \liminf_{j\rightarrow\infty} \frac{
\log\left( (d+1)^{-j} F_j \right)}{j\log d} =
\frac{\log(d+1)-\log \zeta}{\log d}.
\]

This yields the other inequality for $\underline{\alpha}$. The
rest of the proposition follows immediately from Proposition
\ref{prop:alphastarequalsoveralpha}. $\blacksquare$

Another interesting class of examples are the convolutions of
certain Cantor measures. The next lemma is standard but we include the computation for completeness.

\begin{lemma} \label{lemma:convolution}
Let $(p_0,\ldots,p_m)$ and $(p'_0,\ldots,p'_{m'})$ be two
probability vectors. Denote by $\mu$ and $\mu'$ the attractors of
the IFS
\[
\left\{ \left(\frac{x}{d}+\frac{i}{d},p_i\right):0\le i\le m
\right\},\quad \left\{
\left(\frac{x}{d}+\frac{i}{d},p'_i\right):0\le i\le m' \right\}.
\]
(we are not assuming any condition on the weights, $p$, $m$ or
$m'$). Then $\mu*\mu'$ is the attractor of the IFS
\[
\left\{ \left(\frac{x}{d}+\frac{i}{d},v_i\right):0\le i\le m+m',
\right\}
\]
where
\[
v_i=\sum\{p_j p'_{j'}: 0\le j\le m, 0\le j'\le m' \textrm{ and }
j+j'=i\}.
\]

\end{lemma}
\textit{Proof}. Write $\phi_i(x) = x/d + i/d$, $S(x,y)=x+y$, and
compute
\[
\begin{array}{lll}
\mu*\mu'(A)& = & \mu\times\mu'(S^{-1}(A)) \medskip\\
& = & \sum_{i,j} p_i p'_j \mu\times\mu'((S\circ
\phi_i\times\phi_j)^{-1}(A)) \medskip\\
& = & \sum_{i,j} p_i p'_j\mu\times\mu'( (\phi_{i+j}\circ
S)^{-1}(A)) \medskip\\
& = & \sum_{i,j} p_i p'_j (\mu*\mu')\phi_{i+j}^{-1}(A).
\end{array}
\]
The lemma follows. $\blacksquare$

Note that the above lemma also shows that the class of measures
studied here is closed under convolution. We will now briefly
consider the convolutions of biased (middle-third) Cantor
measures. Fix $0<p\le 1/2$, and let $\mu^0_p$ be the attractor of
the IFS $\{(x/3,p),((x+1)/3,1-p)\}$. The lemma shows that if we
let $p_i = {m\choose i} (1-a)^i a^{m-i}$ and $d=3$, then $\mu$ is
the $m$-fold convolution of $\mu^0_p$; note, however, that the
resulting weight vector $\mathbf{p}$ is not always regular. Since
we are assuming $p\le 1/2$, the minimum non-extreme weight is $m
p^{m-1}(1-p)$; it follows that $\mathbf{p}$ is regular if and
only if
\[
(1-p)^m \le m p^{m-1}(1-p) \quad \Longleftrightarrow \quad p \ge
\frac{1}{1+m^\frac{1}{m-1}}.
\]
The numerical value of the minimal $p$ that makes $\mathbf{p}$
regular for $m=3$ is $0.366025$.

Let $\mu_{k}$ be the $k$-fold convolution of the standard Cantor measure. We know from Theorem \ref{theo:extremedimensions} that $\max\Delta\mu_k = k\log 2/\log 3$; because of Proposition \ref{prop:alphastarequalsoveralpha}, this maximum local dimension is isolated. Let

\begin{eqnarray}
\underline{\alpha}_{k} & = & \inf \Delta\mu_{k}; \nonumber\\
\alpha_{k}^{*} & = & \sup \Delta \mu_{k} \backslash \{  k \log 2 / \log 3 \}; \nonumber\\
\gamma_k & = & \dim_H(\mu_k)\nonumber.
\end{eqnarray}

In \cite{convolution} the authors computed the exact values of $\alpha_3^*, \underline{\alpha}_3$ and $\underline{\alpha}_4$. Here we compute the value of $\alpha_4^*$:
\begin{lemma}
\[
\alpha_4^* = \frac{\log(16/5)}{\log 3} = 1.05875...
\]
and this dimension is attained at $x=1/2=\pi(1,1,\ldots)$
\end{lemma}
\textit{Proof}. A calculation shows that $\rho(M(1))=\log(16/5)/\log 3$ and therefore, by Lemma \ref{lemma:periodic}, $\dim\mu_4(1/2) = \log(16/5)/\log 3$.

We claim that
\[
\sigma_1\notin \{0,4\} \,\Longrightarrow\,\eta(\sigma) \ge \frac{1}{5} \left(\frac{5}{16}\right)^{|\sigma|}.
\]
Assuming the claim, the lemma follows at once from the previous calculation and Corollary \ref{coro:approximationalpha}.

The claim will be proved by induction in $|\sigma|$. It is clear for $|\sigma|=1$.  Now let $|\sigma|=j+1>1$ and assume the case $|\sigma|=j$ has been verified. If $\sigma_{j+1}=0$ then $\sigma$ can also be represented by a sequence ending in $3$: since $\sigma_{j+1}=0$ and $\sigma_1>0$, $\pi(\sigma)=i 3^{-j}$ for some positive integer $i$. Then if $\sigma'\in\mathcal{A}^j$ is such that $\pi(\sigma')= (i-1) 3^{-j}$, we get $\pi(\sigma)=\pi(\sigma',3)$. Therefore
\[
\eta(\sigma) \ge (p_0+p_3) \min \left\{ \eta(\sigma'): \sigma'\in\mathcal{A}^j,\sigma_1\notin\{0,4\} \right\}\ge  \frac{5}{16} \frac{1}{5}\left(\frac{5}{16}\right)^{j},
\]
using the inductive hypothesis and the values $p_0 = 1/16; p_3 = 4/16$. A similar argument holds if $\sigma_{j+1} = 1, 3$ or $4$. If $\sigma_{j+1}=2$ then the same conclusion is still true since $p_2 = 6/16 > 5/16$. In any case, the next step is verified and the lemma follows. $\blacksquare$

We remark that although there is no barrier digit for the $4$-fold convolution of the Cantor measure we can still prove that
\[
\Delta\mu_4 = [\underline{\alpha}_4,\alpha_4^*] \cup \{\overline{\alpha}_4\}.
\]
The idea is as follows: let $\sigma^{k,l}$ be the sequence consisting of $k$ twos followed by $l$ ones. By computing the matrices $M(\sigma^{k,l})$ explicitly one can show that the local dimensions arising from such sequences are dense in $[\underline{\alpha}_4,\alpha_4^*]$; by an approximation argument one can show that actually all intermediate dimensions occur.

With the help of Theorem \ref{theo:extremedimensions} (or more precisely Corollary \ref{coro:approximationalpha}) it is possible to estimate the values of $\underline{\alpha}_k$ and $\alpha_k^*$ for other values of $k$, as well as the corresponding parameters in the biased case. It turns out that it is easier to obtain good estimates for $\underline{\alpha}_k$; The following table summarizes these estimates for $5\le k\le 10$. As is to be expected, the smallest local dimensions approach $1$ as $k$ increases, reflecting the progressive smoothing produced by the successive convolutions.
\smallskip

\begin{center}\label{table:mindim}
\begin{tabular}{|c|c|c|}
\hline
$k$ & $\underline{\alpha}_k$ (l.bound) & $\underline{\alpha}_k$ (u.bound) \\
\hline
$5$ & $0.972510$ & $0.972638$ \\
$6$ & $0.976057$ & $0.976628$ \\
$7$ & $0.993697$ & $0.993848$ \\
$8$ & $0.994940$ & $0.995246$ \\
$9$ & $0.998585$ & $0.998657$ \\
$10$ & $0.998908$ & $0.999022$ \\
\hline
\end{tabular}
\end{center}

\smallskip

We remark that in \cite{secondorder} the value of $\gamma_3$ was computed with 4 decimal digits of accuracy (In fact, what they computed is the value of $\tau'(1)$ where $\tau$ is the $L^q$-spectrum of the $3$-fold convolution of the Cantor measure; but it is known that this is equal to the Hausdorff dimension when the measure has an almost sure local dimension, which is the case by Proposition \ref{prop:almostsuredimension}). Theoretically, it is possible to use Proposition \ref{prop:almostsuredimension} to estimate the value of $\gamma_k$, but unfortunately that seems to require extreme computing power (it is not hard to obtain an accuracy of two decimal digits, but since for $k\ge 4$ the values appear to be  between $0.99$ and $1$, this is rather meaningless).

\section{Remarks and open questions}

We finish the paper with some remarks on the relationship between
our results and other recent research in the multifractal theory
of self-similar measures.

\begin{enumerate}
\item
The measures studied here verify the ``weak separation
condition'' (w.s.c.) introduced in \cite{wsc}. The main result of
that paper is that, under the w.s.c., $f(\alpha) = \tau^*(\alpha)$
for any $\alpha=\tau'(q)$, $q\ge 0$. The authors were able to
check the differentiability of $\tau(q)$ in the range $q\ge 0$
for some concrete classes of measures, including a small subset
of the measures analyzed in this paper \cite{secondorder}. Our
examples show that $\tau(q)$ may not be differentiable for $q<0$,
even when $\mu$ is singular. Moreover, this leads to the failure
of the multifractal formalism, as $\tau(q)$ does not ``see'' an
interval of local dimensions. However, it remains a challenging
open question whether $\tau(q)$ is differentiable in $[0,\infty)$
for every weakly-separated self-similar measure.
\item
Olsen \cite{olsen} introduced a more general multifractal
framework. He considers several ``coarse'' and ``fine'' versions
of both the multifractal and $L^q$-spectra; the coarse versions
are defined in terms of coverings or packings by balls of the same
radius, while in the fine spectra variable radii are allowed.
Could it be that the failure of the multifractal formalism in our
setting is due to our consideration of inappropriate
(non-matching) versions of the relevant spectra, rather than to
an intrinsic characteristic of the measures? Even if this were the
case, it would be in sharp contrast with the non-overlapping
situation, in which all of the spectra coincide. We believe,
however, that the breakdown of the multifractal formalism is
essential, and should be corroborated by any reasonably-defined
version of the spectra.
\item
When a barrier digit $b$ is present, the self-similar measure
$\mu$ can also be obtained as the attractor of an
\textit{infinite} iterated function system \textit{without
overlaps}. Indeed, for $\sigma\in\ma^k$ write
\[
\overline{\sigma} = [\pi(\sigma)]_{k}; \quad \Omega = \{
\overline{(\sigma,b)} : \sigma\in (\ma\backslash\{b\})^* \}.
\]
Let
\[
\phi(\overline{\sigma})(x) = \frac{x}{d^{|\sigma|}} + \pi(\sigma);
\]
it is clear that the definition is independent of the
representative chosen. Define also
\[
\mathcal{I} = \left\{ \left(\phi(\overline{\sigma}), \eta(\sigma)
\right) : \overline{\sigma}\in \Omega \right\}.
\]
It is easy to see that $\mu$ is the attractor of the infinite IFS
$\mathcal{I}$ and, since $b$ is a barrier, the open set condition
is verified (take $(0,\xi)$ as the open set). The multifractal
theory for such (and far more general) infinite IFS was developed
in \cite{infiniteifs}, where the multifractal formalism was shown
to hold in certain region, depending on some conditions. Since we
are mainly concerned with the region where the multifractal
formalism fails, we have not attempted to use the results of
\cite{infiniteifs}.

\item The results of this paper give no information about whether
the extreme local dimensions $\underline{\alpha}$ and $\alpha^*$
are attained. We do not know of any example where either of them
is not attained, and we conjecture that this cannot happen.

\item Although our methods do not seem to generalize to other self-similar
measures, our results may still hold in greater generality. In
particular, the following question arises naturally: let $\mu$ be
a self-similar measure on $\mathbb{R}^n$ whose attractor has
nonempty interior. Let $K$ be a compact subset of the interior of
$\supp\mu$. Is it always true that $\mu|_K$ verifies the
multifractal formalism?

We conjecture that the answer is affirmative, at least in the
case where $\supp\mu$ is a linear interval.

\end{enumerate}

\textbf{Acknowledgments}. The author gratefully thanks Doctor
Ursula Molter for many valuable discussions, suggestions and
comments on early version of the paper. He also would like to
thank Professor Boris Solomyak for some helpful conversations and
for pointing to some relevant literature.

\end{document}